\documentclass[11pt,a4paper]{article}
\usepackage[latin1]{inputenc}
\usepackage[german,english]{babel}
\usepackage{amsmath}
\usepackage{amsfonts}
\usepackage{amssymb}
\author{S.J. Patterson}
\title{Rational quadrilaterals and tetrahedra -- the involutions of Kummer, Heegner and Regge}

\date{}

\begin{document}

\maketitle

\section{Introduction}

The purpose of this note is to present some details about involutions
on the set of rational quadrilaterals (quadrilaterals for which the 
distance between every two points is rational) and rational tetrahedra
(set of four points in 3-space for which the distance between any two 
points is rational and moreover the volume is rational).   We can consider 
the case of rational quadrilaterals as tetrahedra of volume $0$.  
One approach to this problem is through the so-called Cayley-Menger 
function, the homogeneous polynomial $CM$ of degree $6$ defined as follows: 
Let $CM_0(x_{12},x_{13},x_{14},x_{23},x_{24},x_{34})$ denote the determinant 
of the matrix
\[
\left[\begin{array}{ccccc}
0&1&1&1&1\\
1&0&x_{12}^2&x_{13}^2&x_{14}^2\\
1&x_{12}^2 &0&x_{23}^2&x_{24}^2\\
1&x_{13}^2 &x_{23}^2&0&x_{34}^2\\
1&x_{14}^2&x_{24}^2&x_{34}^2&0
\end{array}
\right].
\]
This turns out to be a sum of $22$ monomials of degree 6 the coefficients 
of which are $\pm2$; it is therefore convenient to write
\[
CM(x_{12},x_{13},x_{14},x_{23},x_{24},x_{34})=
CM_0(x_{12},x_{13},x_{14},x_{23},x_{24},x_{34})/2 .
\]
If now we have a set of four points, labelled by $1,2,3,4$ and denote 
the distance between $i$ and $j$ by $d_{ij}$  (with $ij$ and $ji$ identified) then, 
if the volume of the tetrahedron be denoted by $V$, we have
\[
(12V)^2 = CM(d_{12},d_{13},d_{14},d_{23},d_{24},d_{34}).
\]
We can regard the two problems as the search for the rational solutions of
$ CM(x_{12},x_{13},x_{14},x_{23},x_{24},x_{34})=0$ (rational quadrilaterals)
and \newline 
$ CM(x_{12},x_{13},x_{14},x_{23},x_{24},x_{34})=y^2$ (rational tetrahedra).
One should be aware that by no means all diophantine solutions correspond to 
geometrical solutions - there are further conditions to be satisfied such as 
the appropriate triangle inequalities and a similar condition for the 
angles at a vertex.  We shall not go into these questions here.

In his paper \cite{EEK} of 1848 Kummer associated to any rational quadrilateral
a degenerate quadrilateral, i.e. one where three vertices lie on a line.   These
can be parametrized rationally.  Kummer then shows that all the rational
quadrilaterals associated with this degenerate quadrilateral are parametrized 
by the rational points on an elliptic curve.   He shows, by a direct construction,
that there are at least two non-trivial rational points on this curve (so it is 
an elliptic curve over $\mathbb Q$ and generally of positive rank).  There is 
a fairly simple constructive method for finding  degenerate quadrilaterals.
It turns out that the associated elliptic curves are often of rank $>1$.  
Once the generators of the group of rational points on the elliptic curve is known
one can construct the corresponding rational quadrilaterals.  This method is 
to a certain degree constructive if one uses the ``generic'' rational points
given by Kummer's construction but it need not be exhaustive.   Furthermore   
there are $12$ ways of carrying out this construction and so it is difficult to 
see a structure in the set.   In 1960 Mordell \cite{M1} returned to Kummer's method 
and showed that generically Kummer's solution is of infinite order and uses this 
to show that rational quadrilaterals are dense in the sense that if a 
quadrilateral is given a rational one can be found arbitrarily close to it.
Kummer was unable to find a method to find rational tetrahedra
\footnote{On this point we quote Schulz' text  \cite[p.19]{OS}:

Hingegen ist das Tetraederproblem als solches wahrscheinlich 
zuerst von Kummer formuliert worden, der den Mitgliedern des 
`` Berliner Mathematischen Seminars''
h\"aufig die Aufgabe gestellt 
hat, Tetraeder mit rationalen Ma\ss zahlen der Kantenl\"angen und des
Volumens aufzusuchen. Abgesehen von mehreren damals berechneten 
Zahlenbeispielen, die sich durch Kleinheit der Kantenl\"angen 
besonders auszeichnen, gelang bereits in einigen der einfacheren 
Spezialf\"alle die vollst\"andige L\"osung. Ferner wurde sp\"ater von
Herrn Schwarz in einem weitaus schwierigeren Spezialfall
eine partikul\"are L\"osung gefunden.

In contrast the problem of rational tetrahedra was probably first 
formulated by Kummer who frequently posed to the members of the 
``Berliner Mathematischen Seminars'' the problem of finding tetrahedra with 
rational side lengths and rational volumes.   Apart from several numerical
examples, characterized by their small side lengths, some progress was made 
in some special families.  Later Herr Schwarz found in a much more 
difficult family.
}

The first method found was due to Friedrich von Ankum
(unpublished) and is analogous to Kummer's method for quadrilaterals\footnote{
For the benefit of the reader we give a brief sketch of how the method functions.
Let $a,b,c$ be the side-lengths of a triangle; we consider the tetrahedroid
$CM(a,b,x,c,y,z)=0$ and suppose that $(x^\circ,y^\circ,z^\circ)$ is such
that $CM(a,b,x^\circ,c,y^\circ,z^\circ)=V^2$. In our case we could take here
a degenerate quadrilateral.   There are $16$ nodes of the tetrahedroid; 
we take one, for example $(0,a,-b)$.   Now consider $CM(a,b,x,c,y,z)$ 
restricted to the line $t \mapsto t(x^\circ,y^\circ,z^\circ)+(1-t)(0,a,-b)$.
As $(0,a,-b)$ is a node this has a double zero at $t=1$; moreover it is 
of degree $4$.   It follows that 
\[CM(a,b,tx^\circ,ty^\circ+(1-t)s,tz^\circ-(1-t)b)=(1-t)^2(R_2t^2+R_1t+R_0)\]
for some polynomials $R_2,R_1,R_0$ in $a,b,c$.  In particular $R_0=V^2$.
This reduces the problem in hand to finding values of $t$ where 
$R_2t^2+R_1t+R_0$.   We know that there is at least one solution and, over
the rational numbers, if the polynomial is neither positive nor negative 
definite it is elementary to find infinitely many solutions.  The exceptional 
case can arise but we have sufficiently many nodes that one generally finds 
a solution.}   
 
 The construction does not, as far as is known, produce all rational
tetrahedra. It leads to a rational curve (instead of an elliptic curve) which is easy to handle.

A second method was found by K. Heegner in his student days.   His idea was to use
cartesian coordinates; one first takes one point to be the origin and a second one
of the form $(d_{12},0)$, or, if one prefers, even of the form $(1,0)$.   The last
two points, in the case of quadrilaterals, turn out to be of the form 
$(x,y\sqrt{\Delta)}$ and  $(x',y'\sqrt{\Delta)}$ where $\Delta$ is a square-free
number such that the area of any triangle with three of these points as vertices is
a multiple of $\sqrt{\Delta}$.   That there is such a number was 
proved by Kummer.   The final two points have to satisfy a number of quadratic
diophantine conditions, namely 
 
\[
\begin{array}{lcc}
x^2+y^2\Delta&=&d_{13}^2\\
x'^2+y'^2\Delta & = & d_{14}^2\\
(x-d_{12})^2+y^2\Delta & = & d_{23}^2\\
(x'-d_{12})^2+y'^2\Delta & = & d_{24}^2\\
(x-x')^2+(y-y')^2\Delta & = & d_{34}^2
\end{array}.
\]

If we wish to treat the case of rational tetrahedra (with rational volumes)
then we can use a three dimensional representation with the third point 
determining the $xy$-plane.  We find that the points can be represented as 
$(0,0,0)$,$(d_{12},0,0)$,$(x,y\sqrt{\Delta},0)$ and 
$(x',y'\sqrt{\Delta},z'\sqrt{\Delta})$ with
\[
\begin{array}{lcc}
x^2+y^2\Delta&=&d_{13}^2\\
x'^2+y'^2\Delta & = & d_{14}^2\\
(x-d_{12})^2+y^2\Delta & = & d_{23}^2\\
(x'-d_{12})^2+y'^2\Delta +z'^2\Delta & = & d_{24}^2\\
(x-x')^2+(y-y')^2\Delta+z'^2\Delta & = & d_{34}^2
\end{array}.
\]

Since they are homogeneous it suffices to solve them in integers and it is not 
difficult to construct algorithms, in particular, for searches in regions bounded
by height.  The use of $CM$ proves to be much more troublesome, especially for 
hand computations. Further applications of this representation allow one to
give systematic proofs of the formul\ae\ used by Kummer.  

Heegner was led by these considerations to study in detail the theory of 
indefinite quadratic forms.    The major fruit of this work is his
Habiliationsthesis, \cite{KHHabil}, on the corresponding Massformel. This 
followed earlier investigations of G. Humbert    \cite{GH1921} ( and was more or less 
simultaneous with the related but more more extensive work  of  Siegel, for example in \cite{CLS1936}) .  
Although his initial hope seems to have been to
use these directly with his representation this never came to anything. 
Heegner became aware (around 1930) of the connection of the 
tetrahedroid with the theory of Kummer surfaces and the Jacobians of 
curves of genus 2  -- it was all formulated differently then -- and
he began to follow a different path.    In the study of reducible jacobians
he began investigations of special examples which he sketched in his last paper
\cite{KH1956}.  Here the Massformel can be used, as is done in \cite[Ch.10]{EGM} 
and \cite[Ch.11]{MR}.   There are hints of what Heegner had in mind on the closing
pages of his paper but they are difficult to interpret.

The algorithms give no insight into structure on the set of solutions in either 
case.  It seems that by the time Heegner became a student in Berlin that the focus
in those circles concerned with Kummer's problem -- mainly around H.A.Schwarz -- 
had become more precise and had in view something more akin to the reduction theory
of quadratic forms or of the Markoff equation. The problem is never quite clearly enunciated 
but Heegner's achievement as a step in the direction of a solution, 
was to show that in the case of rational quadrilaterals there is an action of the
symmetric group $S_6$.   We shall sketch a version of Heegner's proof below.   It
has neither been published or rediscovered in the intervening period.   Heegner's
method also yields an action of $S_5$ in the case of rational tetrahedra.

In recent times a different group, namely the Weyl group of the $D_6$ root
system has attracted a lot of attention.   This is based on an idea of the 
theoretical physicist Tullio Regge  \cite{TR} in a quite different context.
We shall not go into the details here but refer the reader to \cite{DL,DR1,DR2}
where these can be found. This group has $23040$ elements and is an extension 
of $S_6$.    However - it has no direct connection with Heegner's construction.
In fact the two constructions can be combined but in neither case is it known 
how large the group of automorphisms (defined over $\mathbb Q$) of the zero set 
of $CM$ or of $CM-y^2$ is.

The circle of ideas around the Regge involutions is rich and varied and, apart
from geometrical applications includes topics in representation and invariant 
theory.   We shall refer the reader merely to the two papers \cite{DR1,DR2} 
where the reader will find references to the literature. 

Whereas the role of the group $W(D_6)$ is a relatively recent discovery; the
theorem of Heegner is in many  ways similar and is much older - according to
Heegner he found it earlier than 1940 and possibly even in the 1920s.  Heegner
reports that his investigations expanded to such an extent that he was unable to 
publish an account in his life-time. Part of his investigations is contained in the 
unpublished manuscript \cite{KH1}.

Heegner's proof of his theorem is geometrical and is rather sketchy; he promises 
an analytical version in Part 2 but this has not survived, if it were ever 
written. The proof below is presumably akin to his ``analytic'' proof but we
have made use of computer algebra and so it diverges from his.

This paper represents a continuation of \cite{SJP} and we shall take over 
notations and results.   For the convenience of the reader we shall recall
that which is salient, in particular a formula of F. Nei\ss\ and, at a point 
where it is needed, a version of the theory of Weddle surfaces tailored to
our present needs. 

We define the Heron function by

\[ 
H(a,b,c)=(a+b+c)(-a+b+c)(a-b+c)(a+b-c);
\]
Heron's formula gives the area of a triangle with sides $a,b,c$ as 
$\sqrt{H(a,b,c)}/4$. 
Now Neiss' formula (an identity) is:
\[
H(d_{12},d_{13},d_{23})H(d_{12},d_{24},d_{14})=D_{12}^2+
(2 d_{12})^2 CM(d_{12},d_{13},d_{14},d_{23},d_{24},d_{34}) 
\]
where $D_{12}=D_{12}(d_{12},d_{13},d_{14},d_{23},d_{24},d_{34})$ is 
\[
\left|\begin{array}{cccc}
0&1&1&1\\
1&0&d_{12}^2&d_{14}^2\\
1&d_{12}^2&0&d_{24}^2\\
1&d_{13}^2 &d_{23}^2&d_{34}^2
\end{array}
\right|.
\]
We can construct $D_{13},D_{14},D_{23},D_{24},D_{34}$ analogously.  We find
\[
D_{12}D_{34}-D_{13}D_{24}=2(d_{13}^2+d_{24}^2-d_{14}^2-d_{23}^2)CM(d_{12},d_{13},d_{14},d_{23},d_{24},d_{34})
\]
with two analogous formul\ae\ for $D_{12}D_{34}-D_{14}D_{23}$ and
$D_{14}D_{23}-D_{13}D_{24}$.

It is worth pointing out that the argument of  $$D_{12}+2d_{12}\sqrt{CM(d_{12},d_{13},d_{14},d_{23},d_{24},d_{34})}\sqrt{-1}$$
is essentially the dihedral angle along the line between $1$ and $2$, i.e. the angle between the planes containing $1$,$2$,$3$ and
$1$,$2$,$4$.

The volume-preserving property of $W(D_6)$ and Heegner's theorem are of the 
nature where one begins with a number of elementary
transformations, in both cases a copy of $S_4$ with two birational automorphisms
(in Heegner's case) or one additional linear automorphims (``the'' Regge 
automorphism).   Both of these operate on sets of angles, in Heegner's 
case the internal angles of the faces, and in the other case the dihedral 
angles, the angles between faces.  This group also operates linearly on the 
coordinates $\mathbf d$. Despite being similar the two results 
are distinct from one another.   Both have something of an air of finality.
The main point of this note is that this aura in the case of rational 
quadrilaterals and tetrahedra is illusory and much larger groups, presumably
finite, are involved.   The present author has no suggestion what these might be. 

\noindent \textit{ Acknowledgements:} First of all I thank Norbert Schappacher 
for his numerous questions and stimulating comments about Heegner's manuscript 
that forced me to unravel as much of it  as I have achieved.  In particular  his
careful reading of earlier versions of this paper were very helpful. 

In carrying out this work I have made extensive use the \texttt{PARI/gp} 
and of the \texttt{GNU coreutils}.   These have made it possible to 
rework Heegner's arguments and in some cases to add more precision.  It   My
thanks are also due to the teams behind these packages.
is also the intention of the author to provide the reader with the 
tools to carry out further investigations into the these diophantine problems.

\medskip

\begin{center}{\Large{\textbf{First part: The involutions of Kummer, Heegner and Regge}}} \end{center}  

\section{The involutions}

Two of the three involutions are defined by simple formul\ae\ but 
we first note that we can define an action of $S_4$ in both cases simply
by permuting the vertices/variables.   We shall understand the word 
``involution'' to indicate one representative modulo conjugation by this group.
In the cases of the Heegner and Kummer involutions one can arrange the definition
so that one triangle -- typically the $123$-triangle  -- remains fixed.   The
Regge involution does not have this property.   

The construction of the Heegner involution can be given as follows:
\[
\mathbf d \mapsto i_H(\mathbf d)
\]
where, if $\mathbf d =(d_{12},d_{13},d_{14},d_{23},d_{24},d_{34})$ then we
define $i_H$ by 
$$i_H(\mathbf d)=(d_{12},d_{13},d_{12}d_{13}/d_{14},d_{23},
d_{12}d_{34}/d_{14},d_{13}d_{24}/d_{14})$$ 
and an elementary computation shows that 
$$CM(i_H(\mathbf d))=\frac{d_{12}^2d_{13}^2}{d_{14}^4}CM(\mathbf d).$$   The point 
here is that the factor is a square. Note that the triangles with the 
sides $\{i_H(\mathbf d)(13),i_H(\mathbf d)(14), i_H(\mathbf d)(34)\}$ is similar 
to the original triangle $124$, and that the one with sides
$\{i_h\mathbf d(12), i_H\mathbf d(14), i_H\mathbf d(24)\}$ is similar to 
the original $134$.

However the triangle 
$\{i_H(\mathbf{ d})(13),i_H(\mathbf{ d})(34),i_H(\mathbf{ d})(14)\}$  is similar to 
a new triangle with sides $\{d_{12}d_{34},d_{13}d_{24},d_{14}d_{23}\}$.  It is a 
particularly elegant observation of Heegner's that these three quantities
appear in Ptolemy's theorem on quadrilaterals inscribed in a circle.   One
consequence is that if the original quadrilateral is inscribed in a circle
then this triangle is a line segment.   This sheds new light on the treatment 
of such rational quadrilaterals by Brahmegupta and Kummer. 

The Regge involution is given by an even simpler expression, now linear. 
We define $i_R$ by 
\[
\begin{array}{lcl}
i_R{\mathbf d}&=&(d_{12},(-d_{13}+d_{14}+d_{23}+d_{24})/2,\\
&&(d_{13}-d_{14}+d_{23}+d_{24})/2,(d_{13}+d_{14}-d_{23}+d_{24})/2,\\
&& (d_{13}+d_{14}+d_{23}-d_{24})/2,d_{34})).
\end{array}
\]
One verifies that $CM(i_R(\mathbf d))=CM(\mathbf d)$ and also that $i_r^2=id$.

We now come to the earliest but yet most complicated involution.   We 
assume that the quadrilateral is not cyclic, i.e. the vertices do not 
lie on a circle.   The following involution will leave the points 
$1$,$2$,$3$ fixed.  We form the circle $K$ through $1$, $2$ and $4$.   
Then $3$ does not lie on $K$.   We join $3$ and $4$ by a line; this meets 
$K$ in $4$ and another point, say, $5$. This is the image of the $4$ under 
the Kummer involution.   The properties of this map are discussed by Kummer on
\cite[p.17ff.]{EEK}; the construction is illustrated in  Fig. 6 of
\textit{loc.cit.}, Taf.I. Then the Kummer involution, denoted by 
$i_K$ maps $(d_{12},d_{13},d_{14},d_{23},d_{24},d_{25})$ to
\[
\left(d_{12},d_{13},
\frac{D_{13}d_{12}d_{34}}{d_{34}D_{12}},d_{23},
\frac{D_{23}d_{12}d_{14}}{d_{34}D_{12}},
\frac{d_{34}^2D_{12}+d_{24}^2D_{13}+d_{14}^2D_{23}}{d_{34}D_{12}}\right).
\]
In this case the components are not necessarily positive but as the 
polynomial $CM$ involves only squares of the arguments this is not significant.
The formul\ae\ above need a word of explanation.   For the sake of an easier
discussion we rename, for the moment, the points $1,2,3,4$ as $A,B,C,D$ and
the point to be generated as $E$.   The determination of $AE$ and $BE$ uses
the following considerations.   We apply the sine law to $\triangle ADE$ and note 
that by \textit{ Euclid,Elements III.21},   $\angle AED = \angle ADB$. 
It follows that $AD:AE = \sin (\angle ADB):\sin (\angle ADE)$  As 
$\angle ADE = \angle ADC$ or its complement we can replace the last term
by $\sin (\angle ADC)$.   The sines can themselves be expressed in terms of 
the appropriate Heron functions etc..  Finally we have the relations 
noted above between the Heron functions and the $D_{ij}$.  The case of
$d_{25}$ is directly analogous.  Finally we need Ptolemy's Theorem
\textit{Ptolemy, Almagest, 10} and the appropriate formula 
$CE=\pm DE \pm CD$.   It may seem as if there are a number of choices of
signs needed here but a simple algebraic argument shows that one formula
fits all cases.

\section{Heegner's theorem}

Heegner's theorem states that the set of solutions of $ CM(\mathbf d)=0$ has
a group of automorphisms isomorphic to $S_6$ and that this group is 
generated by the Heegner and Kummer automorphisms alongside the elementary 
transformations.   He sketched a proof in \cite{KH1} and promised a fuller
version in part II of that work which has not survived, if it were ever written.
There are a few issues as to the the precise formulation.   We shall sketch a
version of the proof based on Heegner's ideas but where we have used computer
algebra.   Heegner's own argument was based on more classical geometric ideas
and seems to depend on the particular configuration of points.   We should, however,
emphasize that the entire concept is the work of Heegner.
\medskip

\noindent \textbf{Theorem} \textit{The group of automorphisms of} $CM(\mathbf d)=0$
\textit {generated by} $S_4$ \textit{ and the Heegner and Kummer automorphisms 
is isomorphic to} $S_6 \times\{\pm I\}$.
\medskip

The group $\{\pm I\}$ arises from reflections, and is a matter of the formulation.

The proof that the group generated by $S_4$ and the Regge involution is 
$W(D_6)$ is discussed in detail in \cite{DR1,DR2} where there are also references 
to the earlier literature.   Fundamental in this analysis is the role of the 
dihedral angles of the tetrahedron, that is the angles between the faces.
Heegner's proof is also based on the study of angles, but in his case the angles
between the sides.   Since the group ensuing is smaller the calculations are
somewhat easier.

There are $12$ angles between the sides.   In the case of rational quadrilaterals
there are relations arising from each face, and also each vertex.   One can
verify that the resulting set of angles is determined by $5$ of them.   There
is also a transcendental relation (consistency relation) that holds.

In the case of tetrahedra this approach is much less useful and as the Kummer
automorphism is not defined all we need do is to make use of the Heegner 
automorphism and the Ptolemy triangle which leads us to an action of $S_5$. 

There are various ways of making this explicit.   The method used by the  
present author was to denote the vertices by $A$,$B$,$C$,$D$ and to use
the angles $\angle DAB$, $\angle DBC$, $\angle DCA$, $\angle DBA$, 
$\angle DCB$ and $\angle DAC$.   To make this more concrete on can represent 
the points in the complex plane and then we use the vector

\[
\begin{array}{rcl}
\mathbf{\theta}&=&
(\arg((D-A)/(B-A)),\arg((D-B)/(C-B)),\\
&&\arg((D-C)/(A-C)),\arg((D-B)/(A-B)),\\
&&\arg((D-C)/(B-C)),\arg((D-A)/(C-A)).
\end{array}
\]
Although it is not essential we use the choice $-\pi < \arg \le \pi$.
The only condition on $A$,$B$,$C$,$D$ is that the distances between these 
points are those in the rational quadrilateral.   It follows that the 
conjugates have the same property but in this case the arguments are negated.
This is the origin for the doubling noted above.

One verifies that with 
$\mathbf{\theta}=({\theta}_1,{\theta}_2,{\theta}_3,{\theta}_4,{\theta}_5,{\theta}_6)$
 
$${\theta}_1+{\theta}_2+{\theta}_3 \equiv \pi+{\theta}_4+{\theta}_5+{\theta}_6 \pmod {2\pi}$$

and that 
$$\sin({\theta}_1)\sin({\theta}_2)\sin({\theta}_3)=-\sin({\theta}_4)\sin({\theta}_5)\sin({\theta}_6).$$

The latter is the transcendental relation announced above.

We shall use the first of these to work with the first five coordinates.

We now identify number $A$, $B$, $C$ and $D$ as $1$,$2$,$3$,$4$. Then
if we let $(12)$, $(23)$ and $(34)$ be the usual transpositions in $S_4$
they induce linear maps $t_{12},t_{23},t_{34}$.   Likewise the Heegner 
and Kummer involutions are represented by the linear maps, $i_H$, $i_K$.
One finds 

\[
t_{12}=\begin{pmatrix}
 0&0&0&1&0\cr
 1&1&1&-1&-1\cr
 0&0&0&0&1\cr
 1&0&0&0&0\cr
 0&0&1&0&0\cr
\end{pmatrix}
 \]
 \[
t_{23}=\begin{pmatrix}
1&1&1&-1&-1\cr
 0&0&0&0&1\cr
 0&0&0&1&0\cr
 0&0&1&0&0\cr
 0&1&0&0&0\cr
 \end{pmatrix}
 \]
 \[
t_{34}=\begin{pmatrix}
 0&-1&-1&1&1\cr
 0&-1&0&0&0\cr
 -1&-1&0&1&1\cr
 0&-1&0&1&0\cr
 0&-1&0&0&1\cr
 \end{pmatrix}
 \]
 \[
i_H=\begin{pmatrix}
 -1&-1&-1&1&1\cr
 -1&0&0&0&1\cr
 -1&0&0&1&0\cr
 -1&-1&0&1&1\cr
 -1&0&-1&1&1\cr
 \end{pmatrix}
\]
and
\[ 
i_K=\begin{pmatrix}
 0&-1&0&0&1\cr
 -1&0&0&0&1\cr
 0&0&1&0&0\cr
 -1&-1&0&1&1\cr
 0&0&0&0&1\cr  
 \end{pmatrix}.  
\]     
The strategy adopted by the author was to form the list of $32$ elements
$t_{12}^{i_1}t_{23}^{i_2}t_{34}^{i_3}i_H^{i_4}i_K^{i_5}$ with 
$i_1,i_2,i_3,i_4,i_5 \in \{0,1\}$.   It turns out that the set of words of 
length $4$ in this list has $1440$ distinct elements and that $-I$ is amongst 
these.   As $\det{-I}=-1$ we find that the set of elements of determinant $1$ 
has $720$ elements and it is not difficult to see that this is isomorphic to 
$S_6$.   One can, in fact, read off more information but for the present 
purposes this suffices. 

This verification is presumably neither the most efficient nor most elegant but
we shall make do with it for now.   This sketches a proof of Heegner's theorem.

\noindent\textit{Remark}  A variant of this argument, applied to the $6\times 6$
matrix presentation, can be used to confirm that the group generated by $S_4$, ``the'' Regge involution and the $64$ sign-changes has indeed $23040$ elements.
As we are dealing with a Coxeter group it is not difficult to identify the
group as $W(D_6)$.  Indeed, from the list of elements this fact virtually 
presents itself.

\medskip

\begin{center}
\Large{\textbf{Second part: Other involutions}}  
\end{center}
\section{The duality map}
This has also been called the ``switch'' and is particularly simple 
in this case.  Here we hold $d_{12},d_{13},d_{24}$, i.e. the triangle $1 2 3$ fixed  The corresponding Kummer surface is self-dual, also associated with a 
tetrahedroid, and one can identify the duality explicitly - see \cite[Ch. 4]{CF}. 
The map is given by
\[
(d_{14},d_{24},d_{34}) \mapsto 
\left(d_{13} d_{23} \frac{\partial CM}{\partial d_{14}},
d_{12}d_{23} \frac{\partial CM}{\partial d_{24}},
d_{12}d_{13} \frac{\partial CM}{\partial d_{34}}\right)/D
\]
where
\[
D= d_{13} d_{23}d_{14} \frac{\partial CM}{\partial d_{14}}+
d_{12}d_{23}d_{24} \frac{\partial CM}{\partial d_{24}}+
d_{12}d_{13}d_{34} \frac{\partial CM}{\partial d_{34}}.
\]
We have suppressed an irrelevant minus-sign here.   

\section{The Weddle surface and Hutchinson involutions}
The theory of the Weddle surface tailored to the present situation
has been described in \cite{SJP}.  In that article we gave a survey of Cayley's
theory of the symmetroid.  We shall recall the formul\ae\
for the convenience of the reader.   
\newline 
$q_1=(a+b-c)(a-b+c)$, $q_2=-(a+b+c)(-a+b+c)$,\newline $q_3=-(-a+b+c)(a+b-c)$
and $q_4=(a+b+c)(a-b+c)$ 
We let
\[
W_{\mathbf q}(\mathbf v)=
\begin{array}{|cccc|}
v_1^2&v_1&q_1v_1&q_1\\
v_2^2&v_2&q_2v_2&q_2\\
v_3^2&v_3&q_3v_3&q_3\\
v_4^2&v_4&q_4v_4&q_4
\end{array}\ .
\]
The Weddle surface associated with the parameters $a,b,c$ is the zero 
variety of this function.   It is birationally isomorphic to the 
tetrahedroid $CM(a,b,x_1,c,x_2,x_3)=0$, or, more conveniently in
homogeneous form $$x_0^4 CM(a,b,x_1/x_0,c,x_2/x_0,x_3/x_0)=0.$$   
We write $T_{a,b,c}=x_0^4 CM(a,b,x_1/x_0,c,x_2/x_0,x_3/x_0)$.   
For our purposes we need the maps between the tetrahedroid and the Weddle 
surface.   This we now describe.

\[
\begin{array}{cc}
L_1=-(a-b+c)(c x_0+x_2+x_3), & L_1'=(-a+b+c)(-c x_0+x_2+x_3)\\
L_2=(a+b+c)(c x_0+x_2-x_3), & L_2'=(a+b-c)(-c x_0+x_2-x_3)\\
L_3=2(-c x_1-b x_2+a x_3), & L_3'=2(c x_1-b x_2 +a x_3)
\end{array}
\]
\[ 
16c^2T_{a,b,c}=-((L_1L_1'+L_2L_2'-L_3L_3')^2-4 L_1L_1'L_2L_2').
\]
This formula gives another approach to Neiss' identities.

The map from the locus of $T_{a,b,c}=0$ to the Weddle surface can be given as
$(x_0:x_1:x_2:x_3) \mapsto (wv_1:wv_2:wv_3:wv_4)$ with
\[
\begin{array}{l}
wv_1=L_1\cdot( L_1\cdot L'_1-L_2\cdot L'_2-L_3\cdot L'_3)\\
wv_2=L_2\cdot(-L_1L'_1+L_2\cdot L'_2-L_3\cdot L'_3)\\
wv_3=L_3\cdot(-L_1 \cdot L'_1-L_2\cdot L'_2+L_3 \cdot L'_3)\\
wv_4=2\cdot L_1\cdot L_2\cdot L_3. 
\end{array}
\]
Note that because of the intervention of automorphisms neither this map 
nor the next  is unique.  A map in the opposite direction is given by 
$(v_1:v_2:v_3:v_4) \mapsto (S_0:S_1:S_2:S_3)$ where 
\[
\begin{array}{l}
S_0= 2\cdot (v_1\cdot v_2-v_3\cdot v_4)\\
S_1=(a+b+c)\cdot v_1\cdot v_3-(-a+b+c)\cdot v_1\cdot v_4-(a-b+c)\cdot v_2\cdot v_3\\
\hspace{1cm}-(a+b-c)\cdot v_2\cdot v_4\\
S_2=-2\cdot a\cdot v_1\cdot v_2+(a+b+c)\cdot v_1\cdot v_3-(-a+b+c)\cdot v_1\cdot v_4\\
\hspace{1cm}+(a-b+c)\cdot v_2\cdot v_3+(a+b-c)\cdot v_2\cdot v_4-2\cdot a\cdot v_3\cdot v_4\\
S_3=-2\cdot b\cdot v_1\cdot v_2+(a+b+c)\cdot v_1\cdot v_3+(-a+b+c)\cdot v_1\cdot v_4\\
\hspace{1cm}-(a-b+c)\cdot v_2\cdot v_3+(a+b-c)\cdot v_2\cdot  v_4-2\cdot b\cdot v_3\cdot v_4. 
\end{array}
\]

``The'' Hutchinson involution is the map on the Weddle surfae
\[
(v_1:v_2:v_3:v_4)\mapsto (q_1/v_1:q_2/v_2:q_3/v_3:q_4/v_4).
\]
transferred to the tetrahedroid.

\section{Translations by half--periods}
One knows that the abelian variety associated with a tetrahedroid 
is isogenous to the product of two elliptic curves; in our case 
these are defined over $\mathbb Q[\sqrt{-\Delta})$ and are conjugate 
to one another. For a recent account of this see \cite{EJ}.  The 
translations by the half-periods is given by the standard formul\ae\
(as for example in \cite[\S 21.11]{WW}).   In order that the translations
are defined over $\mathbb Q$, and not over the quadratic field, we require
that the half-periods themselves are defined over $\mathbb Q$ and there are
$4$ of these - the other $12$ being irrelevant here.

For our purposes we have preferred to use a version based on the use of
Weierstrass elliptic functions.  There were numerous accounts in the nineteenth
century of ``the parametrization of the Fresnel wave-surface by elliptic 
functions''.   We have found the version given in \cite{DNL} to be particularly
well-suited to our needs.   The two elliptic curves are complex conjugates to 
one another.   We let $a,b,c, H$ be as above.   It turns out that we work in 
$\mathbb Q(\sqrt{-H})$ which is an imaginary quadratic field.   Then the conjugate
curves are given by 
\[
\wp'^2=4(\wp-e_1)(\wp-e_2)(\wp-e_3)
\]
and 
\[
\tilde\wp'^2=4(\tilde \wp-\widetilde {e_1})
(\tilde\wp - \widetilde {e_2})(\tilde\wp -\widetilde {e_3}).
\]

Then $T_{a,b,c}$ is parametrized by 
\[\begin{array}{rcl}
x_0&=&1\\
x_1^2&=&(\wp(u)-e_1)(\tilde\wp(\tilde u)-\widetilde {e_1})\\
x_2^2&=&(\wp(u)-e_2)(\tilde\wp(\tilde u)-\widetilde {e_2})\\
x_3^2&=&(\wp(u)-e_3)(\tilde\wp(\tilde u)-\widetilde {e_3})
\end{array}
\]
for suitable $u,\tilde u$.   There is also a version using the Weierstrass 
$\sigma$-functions which gives the ``theta-function'' formulation.  

We can, using the information in \cite{DNL}, compute the values of the
Weierstrass functions at half-periods in terms of $a,b,c$.   We note at 
the outset that these are only defined up to a scaling factor but that 
this has been decided by the parametrization.  We find 
\[\begin{array}{rcl}
e_1&=&-(a^2-b^2+\sqrt{-H})/3c \\
e_2&=&(a^2-b^2+3c^2+3\sqrt{-H})/6c\\
e_3&=&(a^2-b^2-3c^2+3\sqrt{-H})/6c.
\end{array}
\]
and we can use this to determine a representation of the tetrahedroid as 
a quotient of a product of two elliptic curves.

However there is a much simpler method if we compare the equations above
with Heegner's representation in terms of complex numbers.   Let the
ponits $1,2,3,4$ be represented by $A,B,C,D$; here the triangle  $ABC$
is fixed with side-lengths $a,b,c$.    We are considering $D$ as variable.
We move the origin so that the barycentre of $ABC$ is $0$; i.e. we replace
$A,B,C,D$ by $A^*,B^*,C^*,D^*$ where $A^*=A-(A+B+C)/3$ etc..   Then we can 
take $e_1,e_2,e_3$ to be $A^*,B^*,C^*$ and $\wp(u)=D^*$; $\widetilde{e_1}$ 
etc. to be the complex conjugates of these. 

The action of translation of a half-period $\omega_1$ associated
with $e_1$ is given by the formula
\[\wp(u+\omega_1)=e_1+\frac{(e_1-e_2)(e_1-e_3)}{\wp(u)-e_1}.
\]
This is simply an inversion.   It follows that the three Heegner involutions
associated with $A,B,C$ are the three non-trivial translations associated with 
pairs $(e_1,\bar{e_1})$, $(e_2,\bar{e_2})$ and  $(e_3,\bar{e_3})$.
In these cases the $x_1^2,x_2^2,x_3^2$ will be elements of 
$\mathbb Q(\sqrt{-\Delta})$.  Although these are of interest from the point of
view of algebraic geometry they are irrelevant to our problem.   The show that 
the full group of automorphisms of the tetrahedroid is much larger than the one 
we are considering.

There is now a very extensive theory of the automorphisms of Kummer 
surfaces due to Kond\=o, Keum and Dolgachev. - see \cite{ID1,ID2,KeKo,Ko1,Ko2}.
This theory has also been integrated, by Simon Brandhorst and Matthias Zach into 
the Oscar package for the Julia language.

These automorphisms preserve are of the tetrahedroid, and they are based on 
holding one triangle fixed. We are primarily interested in the automorphisms 
of the higher-dimensional varieties $CM=0$ and of $CM=y^2$.   These have not
been as intensively investigated as those of $K3$ surfaces and the importance
of the Regge transformation is that opens up this area.

\section{Concluding remarks} 

The theory of the group containing the Regge involution, $S_4$ and the 
sign reversals of the coordinates, isomorphic to $W(D_6)$ preserves the
volume of the tetrahedron and it also preserves the rationality of the sides.
The Heegner involution applied to tetrahedra preserves the rationality of the
sides and, given that the sides are rational, multiplies the volume by a 
rational number.   It therefore preserves rationality in the diophantine 
sense but does not belong to $W(D_6)$ which in any case acts linearly on the 
coordinates, whereas the Heegner involution is a birational map.   This means 
that the group generated by ``the'' Heegner involution and the group isomorphic
to $W(D_6)$ is a larger group.   Experiments carried out up to the present 
seem to indicate that the group is finite but do not suggest what this might be.

It will be convenient to gather some facts.    In order to visualise the orbit of
Heegner's group $S_5$ we give representatives modulo the action of $S_4$ 
which is directly geometric.   One such set is (up to similarity)\newline
\[
\begin{array}{l}
(d_{12}d_{14},d_{13}d_{14},d_{12}d_{13},d_{14}d_{23},d_{12}d_{34},d_{13}d_{24}),\\
(d_{12}d_{24},d_{13}d_{24},d_{12}d_{34},d_{23}d_{24},d_{12}d_{23},d_{14}d_{23}),\\
(d_{12}d_{34},d_{13}d_{34},d_{13}d_{24},d_{23}d_{34},d_{14}d_{23},d_{13}d_{23}),\\
(d_{12}d_{34},d_{14}d_{23},d_{14}d_{34},d_{13}d_{24},d_{24}d_{34},d_{14}d_{24})
\end{array}
\] 
along with the original \newline
$(d_{12},d_{13},d_{14},d_{23},d_{24},d_{34})$.
The four faces are the triangles \newline 
$(d_{12},d_{23},d_{13})$
$(d_{12},d_{24},d_{14})$, $(d_{13},d_{34},d_{14})$ and $(d_{23},d_{34},d_{24})$,
again up to similarity.   The faces of the quadrilaterals are a subset of
these combined with the \textit{Ptolemy triangle} 
$(d_{12}d_{34},d_{13}d_{24},d_{14}d_{23})$.   We leave it to the reader to 
tease out the details.  These five triangles, up to similarity, provide 
an invariant under the group generated by $S_4$ and the Heegner involutions.

The same considerations can be applied to the case of rational quadrilaterals
as rational quadrilaterals can be considered as tetrahedra of volume $0$.  
Here it seems as if another finite group, presumably larger, is involved 
in the same way but beyond this indication one cannot at the moment 
even speculate further.

Finally we shall give some illustrative examples to guide the reader.  As 
a basis we shall use the ``smallest'' rational tetrahedron, the one given
as number $19$ is Schulz' table at the end of his dissertation \cite{OS}.
This is 
\[
\mathbf d=(d_{12},d_{13},d_{14},d_{23},d_{24},d_{34})=(2,5,7,6,8,4).
\]
One finds that $CM(\mathbf d)=72^2$ which shows that the volume is $\frac{1}{6}$.
As a companion to this, with the same ``basis'' triangle $(2,5,6)$ (with 
$\Delta = 39$) we may take the rational quadrilateral (found by Heegner's
method) 
\[
\mathbf d = (2,5,20/7,6,8/7,41/7).
\]
Neither of these has symmetries or special features so, although they are 
small, they may be taken as sufficiently ``generic''.   We shall discuss the 
orbits of these under the Heegner and ```Regge'' groups.

We begin with the case of tetrahedra. To give a general picture of the 
situation we need only consider representatives of geometrically similar
(up to a reflection) tetrahedra.   Then the Heegner method provides a 
set of five triangles, the original faces and the ``Ptolemy'' triangle.
Then Heegner's construction is each subset of four triangles yields a 
rational tetrahedron and so we find five examples in all.   For the sake of 
completeness we give such a list, in two forms, one containing the basic triangle 
and the other primitive; apart from $(2,5,7,6,8,4)$ we have 
$(14,35,10,42,8,40)=7(2,5,10/7,6,8/7,40/7)$, $(16,40,8,48,12,42)=
8(2,5,1,6,3/2,21/4)$,\newline
$(8,20,40,24,42,30)=4(2,5,10,6,21/2,15/2)$ \newline 
and the last
one, which does not have $(2,5,6)$ as a face, $(4,21,14,20,16,28)$.   We have
\[
\begin{array}{lcl}
CM(2,5,10/7,6,8/7,40/7)&=&720^2/49^2 \\
CM(2,5,1,6,3/2,21/4)&=&27^2/2^2\\
CM(2,5,10,6,21/2,15/2)&=&135^2 \\
CM(4,21,14,20,16,28)&=&2016^2.
\end{array}
\]

The action of the Regge group is even less spectacular.   The orbit consists of
$23040$ elements but it turns out that they are all doubled.   As the properties 
in which we are interested depend on the variables up to sign we can restrict our 
attention to the non-negative cases and it turns out that there are only 
$216$ of these.   Of these $72$ have at least one entry that is $0$; we 
exclude these as being non-geometric.   Finally we restrict our attention
to those cases where the Heron function of the four faces is positive.  
We find $24$ entries, the orbit under $S_4$ of the tetrahedron with which
we started.

We now turn to the case of rational quadrilaterals.   Heegner's group 
isomorphic to $S_6$ means that we do indeed generate more new and interesting
solutions.  The first level is that of $S_5$ which is simpler than in the 
case of tetrahedra as the parameter $\Delta$ is defined and is the same for 
all faces.  Recall that all the triangles have areas which are rational 
multiples of $\sqrt{\Delta}$; in our case  $\Delta=39$.   The rational
part of the areas is then a useful tool in analysing the faces.

The study of the  $S_5$ action involves the introduction of the Ptolemy triangle
which here is $(20,41,60)$; the Heron function takes the value $99^2\Delta$.
We can consider the $S_5$ action as known; the entire $S_6$ action can 
be summarized by given one representative of each $S_5$ class, that is  giving
just $6=720/120$ rational quadrilaterals.   They can be determined by the method 
of angles and one set, in primitive integral form,  is given by
\[
\begin{array}{l}
(14,35,20,42,8,41),\\
(82,205,128,246,150,99),\\
(164,175,205,66,164,100),\\
(820,875,896,330,132,231),\\
(210,40,70,205,224,,33)\mbox{   and}  \\
(615,375,205,960,656,320).
\end{array}
\] 

The group $W(D_6)$ can be studied as before and now the harvest is marginally
larger.   Apart from the rational quadrilaterals with which we started we find
two new rational quadrilaterals; they are \newline 
$(5/2, 5, 33/14, 11/2, 8/7, 89/14)$
with, again, $\Delta=39$ and
\newline $(5/2, 41/7, 3/2, 65/14, 2, 89/14)$ with now
$\Delta = 1$.   The latter has nothing to do with the  orbit under the
Heegner group of our initial example and is for this reason particularly
interesting. 

In all this we have excluded non-geometric solutions.   It becomes clear
from these examples that this is not the way forward when we study 
the relations between the Heegner and Regge groups.   We leave this to a 
future occasion.

\medskip
\noindent
Mathematisches Institut\\
Bunsenstr. 3--5\\
37073 G\"ottingen\\
Germany

\noindent
e-mail:\texttt{spatter@gwdg.de}

\end{document}